\documentclass[12pt]{amsart}
\usepackage{amsmath,amssymb,amsbsy,amsfonts,amsthm,latexsym,
                        amsopn,amstext,amsxtra,euscript,amscd,mathrsfs,color,bm,cite}
                       
\usepackage{float} 
\usepackage[english]{babel}
\usepackage{mathtools}
\usepackage{todonotes}
\usepackage{url}
\usepackage[colorlinks,linkcolor=blue,anchorcolor=blue,citecolor=blue,backref=page]{hyperref}
\usepackage{enumitem}

\def\leq{\leqslant}
\def\geq{\geqslant}

\newtheorem{theorem}{Theorem}

\newtheorem{lemma}[theorem]{Lemma}

\newcommand{\ord}{\operatorname{ord}}

\newcommand{\Z}{\mathbb{Z}}

\usepackage{etoolbox}

\numberwithin{equation}{section}
\numberwithin{theorem}{section}

\newcommand{\QQ}{\mathbb{Q}}

\def\cS{{\mathcal S}}
\def\cT{{\mathcal T}}
\def\cU{{\mathcal U}}

\def\ov\QQ{\overline{\QQ}}

\def \rad{{\mathrm {rad}}}

\usepackage{mathtools}
\usepackage{todonotes}
\usepackage[norefs,nocites]{refcheck}

\def\({\left(}
\def\){\right)}

\begin{document}

\title[Counting solvable $S$-unit equations]
{Counting solvable $S$-unit equations}

\author[I. E. Shparlinski]{I. E. Shparlinski}
\address{Department of Pure Mathematics, University of New South Wales,
 Sydney, NSW 2052, Australia}
\email{igor.shparlinski@unsw.edu.au}

\author[C. L. Stewart]{C. L. Stewart}
\address{Department of Pure Mathematics, University of Waterloo, 
Waterloo, Ontario, N2L 3G1, Canada}
\email{cstewart@uwaterloo.ca} 
 \begin{abstract}  
We obtain upper bounds on the number of finite sets $\cS$ of primes below a given bound
for which various $2$ variable $\cS$-unit equations have a solution.
\end{abstract}


\maketitle


\section{Introduction}  

Let $s$ be a positive integer and let $\cS=\{p_1,...,p_s\}$ be a set of $s$ distinct primes. The group of $\cS$-units in the rational numbers is the multiplicative group generated by $p_1,...,p_s$ and $-1$. Let $a,b$ and $c$ be non-zero integers. The equation
\begin{equation}\label{unit1}
ax+by=c
\end{equation} 
in $\cS$-units $x$ and $y$ is a two variable $\cS$-unit equation over the rationals. By clearing the denominators from $x$ and $y$ we obtain the equation
\begin{equation}\label{eq:u2w}
au+bv=cw
\end{equation} 
where we ask for solutions in coprime integers $u,v$ and $w$ from the semigroup
$$
\cU_\cS =\{(-1)^{k_0}p_1^{k_1} \ldots p_s^{k_s}:~ k_i = 0,1, \ldots, \ i =1, \ldots s\}.$$
Notice that $\cU_\cS$ is the set of $\cS$-units in the integers.

There is an extensive literature concerning the $\cS$-unit equations~\eqref{unit1} 
and~\eqref{eq:u2w} over the rational numbers and, more generally, over algebraic number fields, see~\cite{EvGy,EGST}. Both~\eqref{unit1} and~\eqref{eq:u2w} have only finitely many solutions, indeed in 1984 Evertse~\cite{Ever} gave an upper bound of  $3 \cdot 7^{2s+1}$ for the number of solutions to~\eqref{eq:u2w}. Further, Erd{\H o}s, Stewart and Tijdeman~\cite{EST} proved that there are arbitrarily large sets $\cS$ for which~\eqref{eq:u2w} with $a=b=c=1$ has at least $\exp\((4+o(1))(s/ \log s)^{1/2}\)$ coprime solutions. Konyagin and Soundararajan~\cite{KoSo} improved the lower bound to $\exp\(s^{2-\sqrt{2}+o(1)}\)$, see also~\cite{Harp,HaSo,LaSo}.

We say that a triple $(a,b,c)$ of non-zero integers is {\it $\cS$-normalized\/} if $a,b,c,p_1,\ldots,p_s$ are pairwise coprime and $0<a\leq b\leq c$.

 Two triples $(a_1,a_2,a_3)$ and $(b_1,b_2,b_3)$ of non-zero integers are said to be 
{\it $\cS$-equivalent\/} if there exists a permutation $\sigma$ of $(1,2,3)$,  a non-zero rational $\lambda$ and $\cS$-units $\varepsilon_1,\varepsilon_2, \varepsilon_3$
such that
$$
b_i=\lambda\varepsilon_ia_{\sigma (i)}
$$
for $i=1,2,3$.

 Each $\cS$-equivalence class in $({\Z}^{*})^3$ contains exactly one $\cS$-normalized triple, see~\cite[p.129]{EGST}. Evertse, Gy{\H o}ry, Stewart and Tijdeman~\cite{EGST} proved that there are only finitely many $\cS$-normalized triples $(a,b,c)$ in  $({\Z}^{*})^3$ for which~\eqref{eq:u2w} has more than two coprime solutions with $w$ positive. The proof depends on Schmidt's Subspace Theorem and as a consequence it does not yield an upper bound on the entries of the normalized triples for which~\eqref{eq:u2w} has more than two coprime solutions with $w$ positive or the size of a solution for such a triple
 
  Put
$$
P= \max\{p_1,\ldots,p_s\}.
$$
 Evertse, Gy{\H o}ry, Stewart and Tijdeman~\cite{EGST} were able to prove an effective result under a more stringent hypothesis. They proved that if $(a,b,c)$ is an $\cS$-normalized triple for which~\eqref{eq:u2w} has at least $s+3$ coprime solutions $(u,v,w)$ with $w>0$ then there exist effectively computable positive numbers $C_1$ and $C_2$ such that  
$$
\max\{a,b,c\} < \exp(s^{C_1s}P^2)
$$ 
and each such solution satisfies
$$
\max\{|u|,|v|,|w|\} < \exp(s^{C_2s}P^3).
$$

Here we adopt a dual point of view.
Instead of fixing $\cS$ and varying $(a,b,c)$ we now fix $(a,b,c)$ and vary $\cS$. We are interested in determining how frequently sets $\cS$ of $s$ primes yield a solution to~\eqref{unit1} and so also~\eqref{eq:u2w}. Similar questions have been raised for other Diophantine equations, see~\cite{AkBha,Brow,BrDi,DiMa} and references therein. However in the context of $\cS$-unit equations this appears to be new. \

 Accordingly, for positive integers $a,b,c,s$ and $H$ we define $N_{a,b,c}(s,H)$ to be the number of $s$ element subsets $\cS=\{p_1, \ldots,p_s\}$ of the primes up to $H$ which are coprime with $abc$ and for which~\eqref{eq:u2w} has a solution in positive integers $u,v,w$ from  $\cU_\cS$.

 Observe that if~\eqref{eq:u2w} has a solution in  $\cU_\cS$ then it also has a solution in $\cU_{\cT}$ for any finite set $\cT$ containing $\cS$. Therefore if $a+b=c$ then $(u,v,w)=(1,1,1)$ is a solution of~\eqref{eq:u2w} and consequently for every finite set of primes $\cS$ the equation~\eqref{eq:u2w} has a solution in positive integers from $\cU_\cS$. For any positive integer $n$ let $\omega(n)$ denote the number of distinct prime factors of $n$ and let $\pi(n)$ denote the number of primes of size at most $n$. Put $\omega(abc)=r$. Thus, provided that $H$ exceeds the greatest prime factor of $abc$,  by the prime number theorem we have
$$
N_{a,b,a+b}(s,H) = \binom{\pi(H)-r}{s} \sim \frac{1}{s!}\(\frac{H}{\log H}\)^s.
$$

If $a=1,b=1,c=1$ and $\cS$=\{2\} then $(u,v,w)=(1,1,2)$ is a solution of~\eqref{eq:u2w} and so every finite set of primes $\cT$ containing the prime $2$ has a solution in positive integers from $\cU_{\cT}$. Further this condition characterizes the finite sets of primes $\cS$ for which~\eqref{eq:u2w} has a solution in positive integers from $\cU_\cS$
 since at least one of $u,v$ and $w$ must be even. Therefore,
$$
N_{1,1,1}(s,H)  = \binom{\pi(H)-1}{s-1} \sim \frac{1}{(s-1)!}\(\frac{H}{\log H}\)^{s-1}.
$$

 We also note that if $a,b$ and $c$ are odd then any set $\cS$ of primes for which~\eqref{eq:u2w} has a solution must include the prime $2$. Therefore in this case 
\begin{equation}\label{eq:PP1}
N_{a,b,c}(s,H) \leq \binom{\pi(H)-1}{s-1} .
\end{equation} 
In general, given a triple $(a,b,c)$ of positive integers we do not know how to characterize the finite sets $\cS$ of primes for which~\eqref{eq:u2w} has a solution. We conjecture however that if $a+b\neq c$ then there is a positive number $C$, which depends on $a,b,c$ and $s$, such that
\begin{equation}\label{eq:PP2}
N_{a,b,c}(s,H) \leq C\(\frac{H}{\log H}\)^{s-1} .
\end{equation} 
Notice that this follows if $a,b$ and $c$ are odd from~\eqref{eq:PP1}  and the prime number theorem.  We  show that in general this follows as a consequence of the {\it $abc$-conjecture\/}, which we state below. We are also able to prove unconditionally an upper bound for $N_{a,b,c}(s,H)$ of the strength of~\eqref{eq:PP2} apart from logarithmic factors if in addition $v$ is required to be small compared to $u$.

Let $\delta$ be a real number with $0\leq \delta \leq 1$ and for positive integers $a,b,c,s$ and $H$ we define $N^{\delta}_{a,b,c}(s,H)$ to be the number of $s$ element subsets $\cS=\{p_1,\ldots,p_s\}$ of the primes up to $H$ which are coprime with $abc$ and for which~\eqref{eq:u2w} has a solution in positive integers $u,v,w$ from $\cU_\cS$ with
\begin{equation}\label{eq:PP3}
v \leq u^{\delta}.
\end{equation} 
Notice that 
\begin{equation}\label{eq:PP4}
N_{a,b,c}(s,H) \leq N^{1}_{a,b,c}(s,H) + N^{1}_{b,a,c}(s,H)
\end{equation}

 For a positive integer $n$ we define $\rad(n)$ to be the greatest squarefree factor of $n$ so
$$
\rad(n) = \prod_{\substack{p \mid n\\p~\text{prime}}} p.
$$
The {\it $abc$-conjecture\/} is that for each positive real number $\varepsilon$ there is a positive number $C(\varepsilon)$ such that if $a,b$ and $c$ are coprime positive integers with $a+b=c$ then
$$
c < C(\varepsilon)\rad(abc)^{1+ \varepsilon},
$$
see, for example,~\cite{StTi, StYu}.

\begin{theorem} \label{thm:Uncond N-eta}
Let $a$, $b$ and $c$ be positive integers with $a+b \neq c$ and let $s$ and $H$ be positive integers with $H \geq 16$. Let $\delta$ be a real number with $0 \leq \delta < 1$. 
\begin{itemize}
\item[(i) ] There is a positive number $C_0$, which is effectively computable in terms of $a$, $b$, $c$, $s$ and $\delta$, such that 
\begin{equation}\label{eq:PP6}
 N^{\delta}_{a,b,c}(s,H) <  C_0 H^{s-1}(\log H)^{s+3}(\log\log H)^2.
\end{equation} 
\item[(ii) ] If $a,b$ and $c$ are odd or if the {\it $abc$-conjecture\/} holds there is a positive number 
$C_1$,   
 which depends on $a$, $b$, $c$ and $s$,  such that 
\begin{equation}\label{eq:PP7}
 N_{a,b,c}(s,H) <   C_1\(\frac{H}{\log H}\)^{s-1}.
\end{equation} 
\end{itemize}
\end{theorem}  

If we require that each prime from $\cS$ divides at least one of $u,v$ and $w$, so that
\begin{equation}\label{eq:PP8}
\omega(uvw)=s,
\end{equation} 
and that $u,v$ and $w$ are pairwise coprime solutions of~\eqref{eq:u2w} then the situation changes. Such solutions we call full rank solutions and they may be viewed as the analogue of full rank solutions in the case of multiplicatively dependent vectors, see~\cite{PSSS, SSS, St}.

Let $a$, $b$ and $c$ be positive integers and let $s$ and $H$ be integers larger than $1$. Define $M_{a,b,c}(s,H)$ to be the number of $s$ element subsets of the primes up to $H$ which are coprime with $abc$ and for which~\eqref{eq:u2w} has a solution in coprime positive integers $u,v,w$ from $\cU_\cS$ for which~\eqref{eq:PP8} holds. 

Furthermore, let $\delta$ be a real number with $0\leq \delta \leq 1$. Define $M^{\delta}_{a,b,c}(s,H)$ as above but with the additional requirement that a solution $(u,v,w)$ satisfies~\eqref{eq:PP3}. In particular,  we have an analogue of~\eqref{eq:PP4}
and we also observe that $ M^{0}_{a,b,c}(s,H)$ corresponds to the equation $au+b=cw$.  
 
Next we  show that the bounds for $M^{\delta}_{a,b,c}(s,H)$ and $M_{a,b,c}(s,H)$ which follow from Theorem~\ref{thm:Uncond N-eta} can be significantly improved.

For any real number $x$ let $\lfloor x \rfloor$ denote the greatest integer less than or equal to $x$.

\begin{theorem} \label{thm:Uncond M}
Let $a$, $b$, $c$, $s$ and $H$ be positive integers with $H \geq 16$. There is a positive number $C_0$, which is effectively computable in terms of $a$, $b$, $c$ and $s$, such that 
\begin{equation}\label{eq:PP9}
 M^{0}_{a,b,c}(s,H) < C_0\(H(\log H)^{s}\log \log H\)^{\lfloor  s/2 \rfloor}.
\end{equation} 
\begin{itemize}
\item[(i) ]Let $\delta$ be a real number with $0< \delta <1$.  There is a positive number $C_1$, which is effectively computable in terms of $a$, $b$, $c$, $s$ and $\delta$, such that 
\begin{equation}\label{eq:PP10}
 M^{\delta}_{a,b,c}(s,H) < C_1 (H(\log H)^{s}\log \log H)^{\lfloor 2s/3 \rfloor}.
\end{equation} 
\item[(ii)] If the {\it $abc$-conjecture\/} holds then there are positive numbers $C_2$ and $C_3$, 
which depend on $a$, $b$, $c$ and $s$, such that
\begin{equation}\label{eq:PP11}
 M^{0}_{a,b,c}(s,H) < C_2\(\frac{H}{\log H}\)^{\lfloor  s/2 \rfloor}
\end{equation} 
and
\begin{equation}\label{eq:PP12}
 M_{a,b,c}(s,H) < C_3\(\frac{H}{\log H}\)^{\lfloor 2s/3 \rfloor}.
\end{equation} 
\end{itemize}  
\end{theorem}

\section{Bounding the exponents}  
For any prime $p$ and non-zero integer $n$ let $\ord_pn$ denote the exact power of $p$ which divides $n$. Let $a,b,c,s$ and $H$ be positive integers and let $\cS=\{p_1,\ldots,p_s\}$ be a set of $s$ prime numbers which are coprime with $abc$ and of size at most $H$.

Here we formulate our two main technical tools which we prove in Sections~\ref{sec:proof L1} and~\ref{sec:proof L2} 
respectively.  

\begin{lemma}
\label{lem: First} Let $\delta$ be a real number with $0\leq \delta<1$. There is a positive number $C_0$,  which is effectively computable in terms of $a,b,c,s$ and $\delta$, such that if $u,v$ and $w$ are pairwise coprime positive integers from $\cU_\cS$ satisfying ~\eqref{eq:u2w} and~\eqref{eq:PP3} then 
$$
\ord_{p_{i}}uvw < C_0(\log H)^{s+1}\log \log H
$$
for i=1,\ldots,s.
\end{lemma} 

The proof of Lemma~\ref{lem: First} depends upon an estimate for linear forms in the logarithms of rational numbers. It is possible to use $p$-adic estimates for linear forms to estimate $\ord_{p_i}uvw$ when $\delta=1$ however the estimate we obtain is too weak in general. We are able to improve Lemma~\ref{lem: First} and treat the case when $\delta=1$ if we assume the {\it $abc$-conjecture\/}.

\begin{lemma}
\label{lem: Second} If the {\it $abc$-conjecture\/} is true then there is a positive number $C_1$,  which depends on $a$, $b$, $c$ and $s$, such that if $u,v$ and $w$ are pairwise coprime positive integers from $\cU_\cS$ satisfying ~\eqref{eq:u2w}  then 
$$
\ord_{p_{i}}uvw < C_1\frac{\log H}{\log p_i}
$$ 
for i=1,\ldots,s.
\end{lemma}  

\section{An estimate for linear forms in the logarithms of rational numbers}  
For the proof of Lemma~\ref{lem: First} we   make use of the following lower bound for a linear form in the logarithms of rational numbers due to Baker and W{\"u}stholz. For any non-zero rational $\alpha$ we have $\alpha=\frac{a}{b}$ with $a$ and $b$ coprime integers and $b$ positive. 
We put $h(\alpha)= \log \(\max\{|a|,|b|\}\)$.

\begin{lemma}
\label{lem: Third} Let $b_1,\ldots,b_n$ be rational integers with absolute value at most $B(\geq3)$. Suppose that $\alpha_1,\ldots,\alpha_n$ are positive rational numbers and put
$$
\Lambda = b_1\log \alpha_1 + \ldots +b_n\log \alpha_n,
$$
where $\log$ denotes the principal branch of the logarithm. If $\Lambda \neq 0$ then there exists an effectively computable positive number $c$ such that
$$
|\Lambda| > \exp\(-(cn)^{2n}\log B \prod_{j=1}^{n} \max\{h(\alpha_j),1\}\).
$$
\end{lemma}  
 This follows from the main theorem of~\cite{BaW}.

\section{Proof of Lemma~\ref{lem: First}}  
\label{sec:proof L1}
Let $\delta$ be a real number with $0\leq \delta <1$. Let $a,b,c,s$ and $H$ be positive integers and let $\cS$ be a set of $s$ primes $p_1,\ldots ,p_s$ which are coprime with $abc$ and of size at most $H$. We may suppose that $a,b$ and $c$ are pairwise coprime. Suppose that~\eqref{eq:u2w} has a solution in coprime positive integers $u,v$ and $w$ from $\cU_\cS$ with $v\leq u^{\delta}$.  Then
$$
u=p_1^{k_1}\ldots p_s^{k_s}, \quad v=p_1^{\ell_1}\ldots p_s^{\ell_s}, \quad  w=p_1^{m_1}\ldots p_s^{m_s}
$$
with $k_1,\ldots ,k_s$, $\ell_1,\ldots ,\ell_s$ and $m_1,\ldots ,m_s$ non-negative integers.

 Let $c_1,c_2,\ldots$ denote positive numbers which are effectively computable in terms of $a$, $b$, $c$, $s$ and $\delta$. Put
$$
 \Lambda=\log \frac{cw}{au},
$$
and note that
$$
 \Lambda=\log (c/a)+ (m_1-k_1)\log p_1 +\ldots + (m_s-k_s)\log p_s .
$$
By~\eqref{eq:u2w}
$$
 \Lambda=\log \(1+\frac{bv}{au}\)
$$
and by~\eqref{eq:PP3}
$$
0< \frac{bv}{au} \leq \frac{b}{au^{1-\delta}}
$$
hence, since $\log(1+t)<t$ for $t>0$,
\begin{equation}\label{eq:PP17}
0< \Lambda < \frac{c_1}{u^{1-\delta}}.
\end{equation} 
On the other hand by Lemma~\ref{lem: Third}, since $h(p_i)\leq \log H$ for $i=1,\ldots,s$,
\begin{equation}\label{eq:PP19}
 \Lambda >  \exp(-c_2\log B (\log H)^s)
\end{equation} 
where
$$
B=\max\{3, |m_1-k_1|,\ldots,|m_s-k_s|\}.
$$
Since $v$ and $u$ are coprime and $p_1,\ldots,p_s$ are coprime with $b$ we see that
for each $i =1, \ldots, k$ either $m_i=0$ or $k_i = 0$ hence
$$
B=\max\{3, m_1,k_1,\ldots,m_s,k_s\}.
$$  
Thus, by~\eqref{eq:PP17} and~\eqref{eq:PP19},
\begin{equation}\label{eq:PP200}
(1-\delta)\log u < c_3 + c_4\log B(\log H)^s.
\end{equation}

Suppose that $B= \max\{k_1,\ldots,k_s\}$.  In that case $u\geq 2^B$ and so, by~\eqref{eq:PP200},
$$
B < c_5 + c_6\log B (\log H)^s
$$
hence
 \begin{equation}\label{eq:PP20}
\max\{k_1,\ldots,k_s\} < c_7(\log H)^{s}\log \log H .
\end{equation} 
On the other hand if $B=\max\{m_1,\ldots,m_s\}$ then $w\geq 2^B$.  Since $v\leq u$ it follows that $u\geq c_{8}2^B$ and we see from~\eqref{eq:PP200} that
 \begin{equation}\label{eq:PP21}
\max\{m_1,\ldots,m_s\} < c_9(\log H)^{s}\log \log H .
\end{equation} 
Finally we note that
$$
v \leq u=p_1^{k_1}\ldots p_s^{k_s} \leq H^{s \max\{k_1,\ldots,k_s)\}}
$$
hence, by~\eqref{eq:PP20} and ~\eqref{eq:PP21},
 \begin{equation}\label{eq:PP22}
\max\{\ell_1,\ldots,\ell_s\} < c_{10}(\log H)^{s+1}\log \log H ,
\end{equation} 
and the result follows from~\eqref{eq:PP20}, ~\eqref{eq:PP21} and~\eqref{eq:PP22}.

\section{Proof of Lemma~\ref{lem: Second}}  
\label{sec:proof L2}

 Let $c_1,c_2,\ldots$ denote positive numbers which depend on $a$, $b$, $c$ and $s$. Suppose that~\eqref{eq:u2w} holds with $u,v$ and $w$ coprime positive integers from $\cU_{\cS}$. We may assume, without loss of generality, that $a,b$ and $c$ are pairwise coprime. Since by assumption $abc$ is coprime with $p_1,\ldots,p_s$ we see that $au,bv$ and $cw$ are pairwise coprime. Thus, by the {\it $abc$-conjecture\/} with $\varepsilon =1$,
$$
cw < c_1Q^2
$$
where $Q$ is the greatest squarefree factor of $aubvcw$. Thus
$$
Q \leq abc\prod_{\substack{p \mid uvw\\p~\text{prime}}} p
$$
hence
$$
w \leq c_2 \prod_{\substack{p \mid uvw\\p~\text{prime}}} p^2.
$$ 
Since
$$
(uvw)^{1/3} < c_3w
$$
we see that
$$
uvw < c_4 \prod_{\substack{p \mid uvw\\p~\text{prime}}} p^6
$$ 
hence
$$
\prod_{\substack{p \mid uvw\\p~\text{prime}}} p^{\ord_puvw} < c_4H^{6s}.
$$
Thus
$$
\ord_{p_i}uvw < c_5\frac{\log H}{\log p_i}
$$
for $i=1,\ldots,s$ as required.

\section{Proof of Theorem~\ref{thm:Uncond N-eta}}  
We first observe that~\eqref{eq:PP7} holds if $a,b$ and $c$ are odd by~\eqref{eq:PP1} and the prime number theorem.

Let $c_1,c_2,\ldots$ be positive numbers which depend on $a$, $b$, $c$, $s$ and $\delta$. Since $a+b\neq c$ we see that $(u,v,w)=(1,1,1)$ is not a solution of~\eqref{eq:u2w} and so each solution of~\eqref{eq:u2w} has $\omega (uvw)\geq 1$.
For each $t$ with $1\leq t\leq s$,  we determine  an upper bound for the number of sets $\cS$ for which~\eqref{eq:u2w} has a solution in coprime positive integers $u,v,w$ from $\cU_\cS$ with $\omega (uvw)=t$.

Accordingly suppose that $\omega(uvw)=1$.  In this case either $u,v$ or $w$ is a power of a prime $p$. We first treat the case when $\omega (u)=1$ and $\omega(vw)=0$.  Then $ap^k+b=c$ with $p$ a prime and $k$ a positive integer. Notice that $(c-b)/a=p^k$ and so the prime $p$ is uniquely determined by $a,b$ and $c$. Once that prime has been fixed we can freely pick the remaining primes in the set $\cS$. We may argue in a similar fashion when $\omega(v)=1$ and $\omega(uw)=0$ and when $\omega(w)=1$ and $\omega(uv)=0$. This gives a total of at most
\begin{equation}\label{eq:1}
3 \binom{\pi(H)-1}{s-1} < c_1\(\frac{H}{\log H}\)^{s-1} 
\end{equation} 
sets $\cS$ for which there is a solution to~\eqref{eq:u2w} with $\omega(uvw)=1$.

If $\omega(uvw)=2$ and one of $u,v$ or $w$ has exactly two prime factors then, as above, the two prime factors are determined by~\eqref{eq:u2w}. Once they have been determined the remaining $s-2$ primes of $\cS$ may be freely chosen and this gives at most
\begin{equation}\label{eq:2}
 c_2\(\frac{H}{\log H}\)^{s-2}.
\end{equation}  
possible sets $\cS$. If $\omega(uvw)=2$ and neither $u,v$ or $w$ has two distinct prime factors then either $\omega(u)=\omega(v)=1$ and $\omega(w)=0$ or $\omega(u)=\omega(w)=1$ and $\omega(v)=0$ or $\omega(v)=\omega(w)=1$ and $\omega(u)=0$. In each of these cases we may suppose that at least one of the primes exceeds $H^{1/3}$ since if two of the primes are less than $H^{1/3}$ and we pick the remainder freely we have at most 
$$
H^{s-2+ 2/3} = H^{s-4/3}
$$
such sets. Once one prime has been determined and its exponent has been fixed there is a corresponding unique prime determined by equation~\eqref{eq:u2w}.

Suppose that~\eqref{eq:PP3} holds. In this case we may apply Lemma~\ref{lem: First} to get an upper bound for the exponent of any prime dividing $uvw$ and consequently an upper bound for the number of possible sets $\cS$. It is 
\begin{equation}\label{eq:4}
c_3\pi(H)^{s-1}(\log H)^{s+1}\log \log H < c_4H^{s-1}(\log H)^2\log \log H.
\end{equation}  
Next we  suppose that the {\it $abc$-conjecture\/} holds. We may assume that one of the primes dividing $uvw$ is larger than $H^{1/3}$. By Lemma~\ref{lem: Second} it occurs to a power at most $c_5$ and so once the largest prime and its exponent have been determined the second prime is fixed by~\eqref{eq:u2w}. The number of sets $\cS$ of primes in this case is at most
$$
c_6\pi(H)^{s-1} < c_7\(\frac{H}{\log H}\)^{s-1}.
$$

Next we suppose that $\omega(uvw)\geq 3$.  Then at least one of $u,v$ and $w$ has at least two distinct primes except in the case when $\omega(u)=\omega(v)=\omega(w)=1$. Suppose we are in the former case. If $\omega(u)\geq 2$ then after fixing the exponents and the primes dividing $v$ and $w$ the primes dividing $u$ are determined. When~\eqref{eq:PP3} holds we may bound the exponents by Lemma~\ref{lem: First}. The number of possible sets $\cS$ in this case is at most
$$
\sum_{j=0}^{s-2} \sum_{r=0}^{j}\binom{\pi(H)}{r} \binom{\pi(H)}{j-r}  (c_8(\log H)^{s+1}\log \log H)^{s-2}.
$$
Treating also the cases when $\omega(v)\geq 2$ and $\omega(w)\geq 2$ we deduce that the number of possible sets $\cS$ is at most
\begin{equation}\label{eq:6}
c_9\pi(H)^{s-2}((\log H)^{s+1}\log \log H)^{s-2} < c_{10}H^{s-2}(\log H)^{(s+1)(s-2)}.
\end{equation}  
On the other hand, if the {\it $abc$-conjecture\/} holds, the number of such sets $\cS$ is, by Lemma~\ref{lem: Second}, at most
$$
c_{11}\pi(H)^{s-2}(\log H)^{s-2}< c_{12}H^{s-2}.
$$

It remains to treat the case when $\omega(u)=\omega(v)=\omega(w)=1$.  We may assume that at most one of the primes dividing $u,v$ and $w$ is smaller than $H^{1/3}$ in size since otherwise, as before, there are at most
$$
\pi(H)^{s-2}\pi(H^{1/3})^2< c_{13}H^{s-{4/3}}.
$$ 
sets $\cS$ with at least two of the primes smaller than $H^{1/3}$. Once the prime and its power has been determined associated with two of $u,v$ and $w$ the prime associated with the third integer is fixed by~\eqref{eq:u2w}. Thus if the {\it $abc$-conjecture\/} holds we have at most
\begin{equation}\label{eq:9}
c_{14}\pi(H)^{s-3}\pi(H)^2< c_{15}\(\frac {H}{\log H}\)^{s-1}.
\end{equation}  
sets $\cS$. On the other hand, by Lemma~\ref{lem: First}, there are at most
\begin{equation}\label{eq:10}
\begin{split}
c_{16}\pi(H)^{s-3}\pi(H)^2&\((\log H)^{s+1}\log \log H\)^2\\
& < c_{17}H^{s-1}(\log H)^{s+3}(\log \log H)^2
\end{split} 
\end{equation}  
sets $\cS$ for which there is a solution of~\eqref{eq:u2w} for which~\eqref{eq:PP3} holds. 

We now see that the bound~\eqref{eq:10} dominates the bounds~\eqref{eq:1}, \eqref{eq:2}, \eqref{eq:4} and~\eqref{eq:6}
and hence we obtain~\eqref{eq:PP6}.
Furthermore, we also see  that the bound~\eqref{eq:9}   dominates the bounds~\eqref{eq:1}, \eqref{eq:2},  \eqref{eq:4} and~\eqref{eq:6}
and hence  we obtain~\eqref{eq:PP7} and our result follows.  

\section{Proof of Theorem~\ref{thm:Uncond M}}  
Note that the result holds when $s=1$ as in the proof of Theorem~\ref{thm:Uncond N-eta} so we may suppose that $s\geq 2$. 

Let $c_1,c_2,\ldots$ be positive numbers which depend on $a$, $b$, $c$, $s$ and $\delta$. If $\delta =0$ we have $v\leq u^0=1$ and so $v=1$ and equation~\eqref{eq:u2w} becomes
\begin{equation}\label{eq:11}
au+b=cw.
\end{equation}  
We are interested in counting sets $\cS$ of primes of size at most $H$ which don't divide $abc$ for which~\eqref{eq:11} has a solution in coprime positive integers $u$ and $w$ from $\cU_{\cS}$ with $\omega (uw)=s$. Notice that in this case either $u$ or $w$ has at least $s/2$ prime factors. 
We  suppose that it is $u$. Put $\omega(w)=r$. In that case $r$ is at most $\lfloor s/2 \rfloor$. There are at most $\binom{\pi(H)}{r}$ choices for the primes dividing $w$ and, by Lemma~\ref{lem: First}, each prime can occur with an exponent of size at most $c_1(\log H)^{s+1}\log \log H$ hence there are at most
\begin{align*}
\sum_{r=0}^{\lfloor s/2 \rfloor}\binom{\pi(H)}{r}& \(c_1(\log H)^{s+1}\log \log H\)^{r}\\
& \qquad \quad  < c_2H^{\lfloor s/2 \rfloor}(\log H)^{s\lfloor s/2 \rfloor}(\log \log H)^{\lfloor s/2 \rfloor}
\end{align*}  
integers $w$ which occur in a solution of~\eqref{eq:11}. Each one determines at most one solution $u$ of~\eqref{eq:11} and the prime factors of $u$ give the remaining primes in $\cS$. The same argument applies when $u$ has at most $\lfloor s/2 \rfloor$ prime factors. Therefore~\eqref{eq:PP9} follows.

If there is a solution of~\eqref{eq:u2w} with $\omega (uvw)=s$ then at least one of $u,v$ and $w$ has at least $s/3$ distinct prime factors, say $u$ for instance. Then the number of distinct prime factors of $vw$ is at most $\lfloor 2s/3 \rfloor$. The exponents associated with these primes are of size at most $c_3(\log H)^{s+1}\log \log H$ if $0\leq\delta < 1$ by Lemma~\ref{lem: First}. The remaining primes dividing $u$ are determined by~\eqref{eq:u2w} once we know the primes and their exponents dividing $v$ and $w$. The number of possible sets $\cS$ in this case is at most
$$
 \sum_{j=0}^{\lfloor 2s/3 \rfloor}\displaystyle\sum_{r=0}^{j}\binom{\pi(H)}{r} \binom{\pi(H)}{j-r} (c_3(\log H)^{s+1}\log \log H)^{\lfloor 2s/3 \rfloor}.
$$
Arguing as above with $v$ and $w$ in place of $u$ we find that the number of sets $\cS$ is at most
 \begin{equation}\label{eq:13}
c_4\pi(H)^{\lfloor 2s/3 \rfloor}((\log H)^{s+1}\log \log H)^{\lfloor  2s/3 \rfloor}
\end{equation}  
and~\eqref{eq:PP10} now follows from~\eqref{eq:13} and the prime number theorem.

The number of subsets of at most $\lfloor 2s/3 \rfloor$ primes up to $H$ with at least one prime of size at most $H^{1/3}$ is at most
$$
c_5H^{\lfloor 2s/3 \rfloor -1}H^{1/3}
$$
and the number of subsets of at most $\lfloor s/2 \rfloor$ primes up to $H$ with at least one prime of size at most $H^{1/3}$ is at most
$$
c_6H^{\lfloor s/2 \rfloor -1}H^{1/3}
$$  
and, subject to the {\it $abc$-conjecture\/}, the number of possible exponents for each prime is at most $c_7\log H$ by Lemma~\ref{lem: Second}. This then determines at most
$$
c_8H^{\lfloor 2s/3 \rfloor -2/3}\log H
$$
sets in the first case and at most
$$
c_9H^{\lfloor s/2 \rfloor -2/3}\log H
$$
sets $\cS$ in the second case. On the other hand if all of the primes are at least $H^{1/3}$ then, subject to the {\it $abc$-conjecture\/}, the number of possible exponents for each prime is at most $c_{10}$ hence the total number of sets $\cS$ determined is at most
$$
c_{11}\(\frac{H}{\log H}\)^{\lfloor 2s/3 \rfloor }
$$
in the first case so~\eqref{eq:PP12} follows and is at most
$$
c_{12}\(\frac{H}{\log H}\)^{\lfloor s/2 \rfloor }
$$
in the second case so~\eqref{eq:PP11} follows.

\section*{Acknowledgements}

This work started  during a very enjoyable visit of the first author 
to the University of Waterloo, whose hospitality and support are very much appreciated. 

The first author was supported by the Australian
Research Council Grant DP170100786. 
The research of the second author was supported in part by the Canada Research Chairs Program and by Grant A3528 from the Natural Sciences and Engineering Research Council of Canada.

\end{document}